\documentclass[11pt,a4paper]{amsart}
\pdfoutput=1

\usepackage[utf8]{inputenc}

\usepackage{hyperref,graphicx,amssymb}

\DeclareMathOperator{\AutHS}{Aut(HS)}
\DeclareMathOperator{\HS}{HS}
\DeclareMathOperator{\Gal}{Gal}
\DeclareMathOperator{\Sym}{Sym}

\newcommand{\QQ}{\mathbb{Q}}
\newcommand{\CC}{\mathbb{C}}

\newtheorem*{remark*}{Remark}

\begin{document}

\title[Higman-Sims Group as a Galois Group]{Explicit Polynomials Having the Higman-Sims Group as Galois Group over $\mathbb{Q}(t)$}

\author{Dominik Barth}
\author{Andreas Wenz}

\address{Institute of Mathematics\\ University of Würzburg \\ Emil-Fischer-Straße 30 \\ 97074 Würzburg, Germany}
\email{dominik.barth@mathematik.uni-wuerzburg.de}
\email{andreas.wenz@mathematik.uni-wuerzburg.de}



\begin{abstract}
We compute explicit polynomials having the sporadic Higman-Sims group $\mathrm{HS}$ and its automorphism group $\mathrm{Aut(HS)}$ as Galois groups over the rational function field $\mathbb{Q}(t)$.
\end{abstract}

\maketitle

\section{Introduction}

From a theoretical perspective it is known that $\AutHS$, the automorphism group of the sporadic Higman-Sims group $\HS$, occurs as a Galois group over $\QQ(t)$ since it has a rigid rational generating triple, see \cite{Hunt} and \cite{Malle}.

In order to find explicit polynomials with Galois group $\AutHS$ over $\QQ(t)$ one can compute a three-point branched covering, also called \emph{Belyi map}, over $\mathbb{P}^1\CC$ whose ramification corresponds to these rigid rational triples.

For a thorough survey on computing Belyi maps refer to \cite{Sijsling}. Recently, Klug et al.\ calculated a Belyi map of degree 50 with monodromy group isomorphic to $\textnormal{PSU}_3(\mathbb{F}_5)$ using modular forms, see \cite{Klug}.

We developed another efficient method of computing certain Beyli maps of higher degree which we will explain in detail in an upcoming paper.  
The purpose of the current note is to present a Belyi map of degree $100$ with monodromy group isomorphic to $\AutHS$.
As a consequence, we obtain polynomials having $\HS$ and $\AutHS$ as Galois groups over $\QQ(t)$.

\section{Ramification Data and Computed Results}

Our goal is to compute a Belyi map $f: \mathbb{P}^1\CC \rightarrow \mathbb{P}^1\CC$ of ramification type $(x,y,z) \in S_{100}^3$ given by
\begin{align*}
x \;=\; & (1, 64, 8, 54, 37)(2, 20, 81, 42, 49)(3, 98, 32, 73, 89)(4, 96, 86, 15, 79)\\
   & (5, 22, 28, 78, 48)(6, 67, 97, 40, 14)(7, 58, 82, 59, 18)(9, 16, 87, 85, 60)\\
   & (10, 70, 41, 56, 55)(11, 77, 36, 25, 68)(12, 17, 19, 21, 80)(13, 35, 90, 33, 91)\\
   & (23, 50, 66, 84, 27)(24, 72, 95, 52, 76)(26, 99, 100, 57, 93)(29, 71, 38, 69, 65)\\
   & (30, 74, 94, 53, 51)(31, 45, 47, 75, 34)(43, 63, 44, 46, 62),
   \end{align*}
   \begin{align*}
y \;=\; & (1, 20)(2, 64)(3, 76)(4, 45)(5, 83)(6, 26)(7, 13)(8, 74)(9, 41)(10, 63)(11, 25) \\
   &(12, 66)(14, 21)(15, 52)(16, 62)(17, 33)(18, 35)(19, 42)(22, 60)(23, 58)\\    
   &(24, 73)(28, 98)(29, 82)(30, 53)(31, 61)(32, 59)(34, 67)(36, 95)(37, 85)\\
   &(38, 47)(39, 51)(40, 80)(43, 92)(44, 78)(46, 99)(48, 55)(49, 94)(50, 91)\\
   &(54, 90)(65, 88)(69, 72)(71, 75)(77, 79)(81, 87)(84, 97)(86, 100)(93, 96),
\\[2mm]
z \;=\; & (1, 2) (13, 18)(21, 40) (47, 71) (25, 68) (15, 95, 77)(20, 37, 87)(39, 53, 51) \\
 & (49, 74, 64)(55, 78, 63)(57, 100, 96)(66, 80, 97)(73, 76, 89)(7, 91, 23) \\
 & (12, 50, 33) (22, 85, 54, 35, 59, 98)(24, 32, 82, 65, 88, 69)\\
 & (3, 52, 86, 99, 44, 28)(4, 31, 61, 34, 6, 93)(5, 83, 48, 56, 41, 60)\\
&(8, 30, 94, 42, 17, 90)(9, 70, 10, 43, 92, 62)(11, 36, 72, 38, 45, 79)\\
&(14, 19, 81, 16, 46, 26)(27, 84, 67, 75, 29, 58).
\end{align*}
This permutation triple is of the following type:
\renewcommand{\arraystretch}{1.5}
\begin{center}
\begin{tabular}{c|c|c|c}
& $x$ & $y$ & $z$ \\  \hline
cycle structure & $5^{19}.1^5$ & $2^{47}.1^6$ & $6^{10}.3^{10}.2^5$ \\
\end{tabular}
\end{center}
With the help of \texttt{Magma} \cite{Magma} we can easily verify:
\begin{itemize}
\item $x\cdot y \cdot z = 1$
\item $\AutHS =  \langle x,y \rangle$
\item $(x,y,z)$ is a rigid and rational triple of genus $0$
\end{itemize}
Due to the rational rigidity criterion \cite[p.\,48]{voelklein} and a rationality consideration there exists a Belyi map $f\in \mathbb{Q}(X)$ of degree $100$ with monodromy group isomorphic to $\AutHS$. Note that $f$ is unique, up to inner and outer Möbius transformations.

Applying our newly developed method we were able compute this Belyi map explicitly. 
The resulting function $f: \mathbb{P}^1\CC \rightarrow \mathbb{P}^1\CC$ is of the form
\[
f(X) = \frac{p(X)}{q(X)} = 1 + \frac{r(X)}{q(X)}
\]
where
\begin{align*}
p(X) \;=\;3^3 \;\cdot\; & (X^4 - 8 X^3 - 6 X^2 + 8X + 1)^5 \cdot \\
    &(X^5 - 5X^4 + 50X^3 + 70X^2 + 25X + 3)^5 \cdot \\
    &(3X^5 - 5X^4 - 5X^3 + 35X^2 + 40X + 4) \cdot \\
    &(9X^{10} - 30X^9 + 55X^8 - 200X^7 + 210X^6 + 924X^5 \\
    &- 890X^4 - 360X^3 + 1925X^2 - 1070X + 291)^5,
\end{align*}
\begin{align*}
q(X) \;=\;
    &(3X^5 - 35X^4 + 90X^3 - 50X^2 + 15X + 9)^2 \cdot \\
    &(9X^{10} - 120X^9 + 10X^8 - 1960X^7 - 1090X^6 + 3304X^5 \\
    &  - 760X^4 - 920X^3 + 145X^2 + 80X + 6)^3 \cdot \\
    & (3X^{10} - 10X^9 - 65X^8 + 160X^7 - 90X^6 - 932X^5 \\
    & - 330X^4 + 880X^3 + 1255X^2 + 830X + 27)^6
\end{align*}
and
\begin{align*}
r(X) &\;=\; p(X)-q(X)\\ &\;=\; 2^2 \cdot 3^{14} \cdot 5^3 \cdot (X-1)\cdot r_5(X) \cdot r_{10}^2(X) \cdot r_{16}^2(X) \cdot r_{20}^2(X)
\end{align*}
with irreducible monic polynomials $r_j$  of degree $j$.

From the factorizations of $p$, $q$, $r$ and the Riemann-Hurwitz formula it is clear that $f$ is indeed a three-point branched cover of $\mathbb{P}^1\CC$, ramified over $0,1$ and $\infty$.

\section{Verification of Monodromy}

We will present two proofs to verify that the monodromy group of our Belyi map $f=p/q$ is isomorphic to $\AutHS$.

\begin{figure}[htbp]
    \centering
        \includegraphics[clip, trim=2cm 1.2cm 2cm 0.5cm, width=1.00\textwidth]{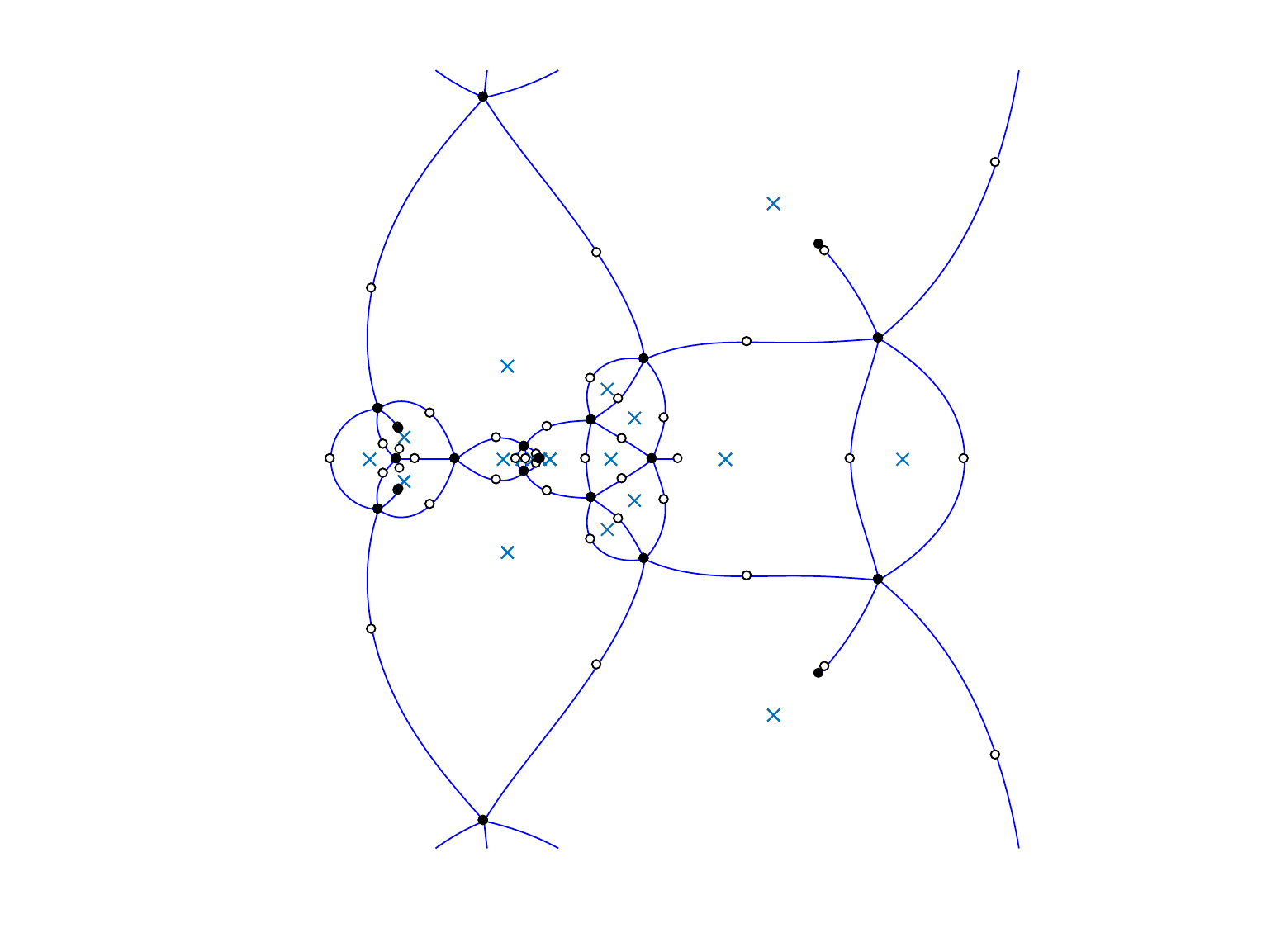}
    \caption{Dessin d'enfant corresponding to $f$}
    \label{fig:dessin}
\end{figure}

First, one can compute the corresponding dessin d'enfant, i.e.\ the bipartite graph drawn on the Riemann sphere $\mathbb{P}^1\CC$ obtained by taking the elements of $f^{-1}(0)$ as black vertices, those of $f^{-1}(1)$ as white vertices and the connected components of $f^{-1}((0,1))$ as edges, labelled from $1$ to $100$. A part of this bipartite graph is shown in Figure~\ref{fig:dessin}. Note that the poles of $f$ are marked by '$\times$'.
Listing the cyclic arrangement of adjacent edges around each black and white vertex, respectively, we obtain the cycles of $x$ and $y$, up to simultaneous conjugation. Thus the monodromy group of $f$ is isomorphic to $\AutHS$.

Another way to verify the monodromy can be done algebraically:
The monodromy group of $f$ can be viewed as the Galois group $\Gal(p(X)-tq(X)\mid \QQ(t))$ or equivalently $\Gal(p(X)-f(t)q(X)\mid \QQ(f(t)))$.

First note that  $\Gal(p(X)-f(t)q(X)\mid \QQ(t))$ equals the point stabilizer of $t$ in the permutation group $\Gal(p(X)-f(t)q(X)\mid \QQ(f(t)))$ acting transitively on the $100$ roots of $p(X)-f(t)q(X)$.

As $p(X)-f(t)q(X)$ factorizes over $\QQ(t)[X]$ into three irreducible polynomials of degrees $1$, $22$ and $77$, respectively, we see that $\Gal(p(X)-f(t)q(X)\mid \QQ(f(t)))$ and thus also $G:=\Gal(p(X)-tq(X)\mid \QQ(t))$ are rank 3 permutation groups of degree $100$ with subdegrees $1$, $22$ and $77$.

We now show that $G$ is actually a primitive permutation group:
Suppose $G \leq \Sym(\Omega)$, $|\Omega|=100$, has some non-trivial block $\Delta$, i.e.\ $1 < |\Delta| < 100$, such that for each $g \in G$ either $\Delta^g = \Delta$ or $\Delta ^g \cap \Delta = \emptyset$. Now fix some $\omega \in \Delta$. Then the stabilizer $G_\omega$ must leave $\Delta$ invariant, and --- as $G$ is rank $3$ group --- $G_\omega$ has exactly the non-empty orbits $\lbrace \omega \rbrace$, $\Delta \setminus \lbrace \omega \rbrace$ and $\Omega \setminus \Delta$.
Knowing the sizes of the suborbits we find that $\Delta$ has either length $1+22=23$ or length $1+77=78$.
This is a contradiction as the size of a block must always divide the permutation degree, in our case $100$.

Now, combining the classification of all finite primitive rank $3$ permutation groups (see e.g.\ \cite{Liebeck}) with the subdegrees of $G$, only two possibilities remain: $G=\AutHS$ or $G = \HS$.

Since $\HS$, in contrary to $\AutHS$, is an even permutation group, it suffices to check whether the discriminant $\delta$ of $p(X)-tq(X) \in \QQ(t)[X]$ is a square in $\QQ(t)$.
Using \texttt{Magma} we see $\delta = u^2 2(t-1)$ for some $u \in \QQ(t)$ and therefore $G=\AutHS$.

\begin{remark*}
By applying the previous arguments to $p(X)-(2t^2+1)q(X)$ we find $\Gal(p(X)-(2t^2+1)q(X)\mid \QQ(t))$ is either $\HS$ or $\AutHS$. The discriminant, however, is a square now, thus $\Gal(p(X)-(2t^2+1)q(X)\mid \QQ(t)) = \HS$.
\end{remark*}

\section{Another Example}

Essentially, $\AutHS$ contains exactly two rigid rational generating triples of genus 0. The first one has been discussed in the previous section. The second triple $(x,y,z)\in S_{100}$ where
\begin{align*}
x \;=\;& (1, 23, 53, 86)(2, 36, 29, 43)(3, 15, 46, 6)(4, 80, 71, 81)(5, 75, 16, 47)
 \\& (7, 32, 60, 8)(9, 76, 100, 51)(10, 50, 49,
    34)(11, 28, 74, 84)(12, 72, 37, 52)
    \\& (13, 21, 96, 88)(14, 41, 40, 87)(17, 42, 45, 79)(18, 63, 19, 20)(22, 99, 39, 
    89)
    \\& (24, 59, 77, 38)(25, 68, 26, 35)(27, 69, 73, 48)(30, 92, 33, 82)(31, 56, 93, 58) \\&(44, 98, 67, 64)(54, 95, 85, 
    62)(55, 65, 94, 61)(57, 78, 83, 97)(66, 90, 70, 91),
\\[2mm]
y \;=\;&(1, 75, 5, 71, 15)(2, 43, 52, 89, 39)(3, 18, 100, 33, 35, 26, 58, 32, 53, 23)
\\ &(4, 81, 47, 16, 86, 7, 42, 38, 77, 
    59)(6, 41, 14, 87, 82, 76, 9, 97, 19, 63)
    \\&(8, 60, 93, 56, 13, 61, 36, 99, 70, 45)(10, 65, 55, 88, 12, 29, 94, 34,
    49, 50)
    \\& (11, 44, 64, 25, 92)(17, 72, 96, 69, 28, 30, 40, 46, 80, 24)(20, 83, 78, 57, 51)\\&(21, 31, 68, 67, 98, 84, 
    74, 27, 48, 73)(22, 37, 79, 90, 66, 95, 54, 62, 85, 91)
\end{align*}
and $z:= (xy)^{-1}$

\noindent of ramification type
\renewcommand{\arraystretch}{1.5}
\begin{center}
\begin{tabular}{c|c|c|c}
& $x$ & $y$ & $z$ \\  \hline
cycle structure & $4^{25}$ & $10^{8}.5^4$ & $2^{35}.1^{30}$ \\
\end{tabular}
\end{center}
leads to the Belyi map
\[
f(X) = \frac{p(X)}{q(X)} = 1 + \frac{r(X)}{q(X)}
\]
where
\begin{align*}
p(X) \;=\;
     &(7X^5 - 30X^4 + 30X^3 + 40X^2 - 95X + 50) ^4 \cdot \\ 
     &(2X^{10} - 20X^9 + 90X^8 - 240X^7 + 435X^6 - 550X^5 \\ &+ 425X^4 - 100X^3 - 175X^2 + 250X - 125)^4 \cdot \\
    & (2X^{10} + 5X^8 - 40X^6 + 50X^4 - 50X^2 + 125)^4,\\
q(X) \;=\;&
2^8 \cdot  (X^4 - 5)^5 \cdot \\
& (X^8 - 20X^6 + 60X^5 - 70X^4 + 100X^2 - 100X + 25)^{10}.    
\end{align*}
Of course, it remains to verify that this rational function is indeed a three-point branched cover having the desired monodromy group. However, this can be done in the exact same way we already demonstrated in the previous section.

\section*{Acknowledgements}

We would like to thank Peter Müller for introducing us to the subject of this work as well as for sharing some tricks to verify the monodromy algebraically,
and Joachim König for providing us with the rigid rational permutation triples of $\AutHS$ we realized as a Galois group over $\QQ(t)$ in this paper.


\begin{thebibliography}{1}

\bibitem{Magma}
Wieb Bosma, John Cannon, and Catherine Playoust.
\newblock The {M}agma algebra system. {I}. {T}he user language.
\newblock {\em J. Symbolic Comput.}, 24(3-4):235--265, 1997.
\newblock Computational algebra and number theory (London, 1993).

\bibitem{Hunt}
David~C. {Hunt}.
\newblock {Rational rigidity and the sporadic groups.}
\newblock {\em {J. Algebra}}, 99:577--592, 1986.

\bibitem{Klug}
Michael Klug, Michael Musty, Sam Schiavone, and John Voight.
\newblock Numerical calculation of three-point branched covers of the
  projective line.
\newblock {\em LMS Journal of Computation and Mathematics}, 17(1):379--430, 001
  2014.

\bibitem{Liebeck}
Martin~W. {Liebeck} and Jan {Saxl}.
\newblock {The finite primitive permutation groups of rank three.}
\newblock {\em {Bull. Lond. Math. Soc.}}, 18:165--172, 1986.

\bibitem{Malle}
Gunter {Malle} and B. Heinrich {Matzat}.
\newblock {\em {Inverse Galois theory.}}
\newblock Berlin: Springer, 1999.

\bibitem{Sijsling}
Jeroen~{Sijsling} and John~{Voight}.
\newblock {On computing Belyi maps.}
\newblock In {\em {Num\'ero consacr\'e au trimestre ``M\'ethodes
  arithm\'etiques et applications'', automne 2013}}, pages 73--131.
  Besan\c{c}on: Presses Universitaires de Franche-Comt\'e, 2014.

\bibitem{voelklein}
Helmut {V\"olklein}.
\newblock {\em {Groups as Galois groups: an introduction.}}
\newblock Cambridge: Cambridge Univ. Press, 1996.

\end{thebibliography}
\end{document}